\documentstyle[aaai]{article}

\newtheorem{Def}{Definition}

\newtheorem{Thm}{Theorem}

\newtheorem{Rem}{Remark}
\def\(-{(\hspace{-1.5 mm}}
\def\s~{| \hspace{-1.8 mm} \sim}
\def\A{{\cal A}}

\def\F{{\cal F}}
\def\I{{\cal I}}
\def\K{{\cal K}}

\def\P{{\cal P}}

\def\R{\cal R}
\begin{document}
\title{Precise Propagation of Upper and Lower Probability
Bounds in System P}
\author{Angelo Gilio \\ Dipartimento di Matematica e Informatica, Citt\`a Universitaria \\
Viale A. Doria, 6 - 95125 Catania (Italy).\\
e-mail: gilio@dmi.unict.it}
\maketitle

\begin{abstract}
\begin{quote}
In this paper we consider the inference rules of System P in the framework of
coherent imprecise probabilistic assessments. Exploiting our algorithms, we
propagate the lower and upper probability bounds associated with the
conditional assertions of a given knowledge base, automatically obtaining the
precise probability bounds for the derived conclusions of the inference rules.
This allows a more flexible and realistic use of System P in
default reasoning and provides an exact illustration of the degradation of
the inference rules when interpreted in probabilistic terms. We also examine
the {\em disjunctive Weak Rational Monotony} of System P$^+$ proposed by Adams
in his extended probability logic.
\end{quote}
\end{abstract}

\ \\

\section{Introduction}
In the applications of intelligent systems to  automated  uncertain
reasoning the explicit knowledge of the agent is represented by  a
knowledge base $\cal K$, constituted by  a  set  of  {\em  conditional
assertions} (i.e. {\em defaults}).  The
nonmonotonic  inferential process is developed  using  a  suitable
set of rules. Among the many formalisms which have been proposed for
default reasoning, the so-called System P developed in
(Kraus, Lehmann, and  Magidor 1990) is
widely accepted and has a probabilistic semantics based on infinitesimal
probabilities, see (Adams 1975, Pearl 1988). An extended
probability logic has been proposed in (Adams 1986) by allowing
{\em disjunction of conditionals}. The corresponding system P$^+$ has been
studied in (Schurz 1998) where the perspectives of the infinitesimal probability
semantics and that of a noninfinitesimal one, based on probabilistic
inequality relations, have been unified.
In (Hawthorne 1996) many logics of nonmonotonic conditionals, that behave like
conditional probabilities at various levels of probabilistic support, have
been examined. In the quoted paper the author shows that, for each
given conditional $\rightarrow$, there is a probabilistic support level $r$
and a conditional probability $P$ such that, for all sentences $B, A$,
it is $B \rightarrow A$ only if $P(A|B) \geq r$.
We recall that an early examination of Adams rules by means of imprecise
probabilities has been given in (Dubois, and Prade 1991). In the quoted
paper a semantic interpretation in terms of intervals has been given for the
relations of negligibility, closeness and comparability examined in the
system of qualitative reasoning proposed in (Raiman 1986). Moreover,
an application to the inference rules of Adams has been given interpreting
"$P(A|B)$ is large" by means of the relation of closeness to 1
(making the infinitesimal parameter $\epsilon$ explicit).
In this way a quantification of
the degradation of the validity of Adams' rules when reasoning with
noninfinitesimal probabilities has been obtained.
While in practical applications the semantics of infinitesimal probabilities
may involve some difficulties, the approach based
on lower (and possibly upper) probability bounds, proposed also in
(Bourne, and Parsons 1998), is clearly more flexible and realistic.
In this paper the propagation of probability bounds to the conditional
assertions associated with the rules of System P is examined
in the framework of the de Finetti's probabilistic
methodology, based on the coherence principle. Notice that
a coherent set of conditional probability assessments satisfies all
the usual properties of conditional probabilities.
A short examination of the logic of conditionals of Adams
from the point of view of coherence has been given in (Gilio 1997),
where the propagation of probabilistic bounds has not been considered.
We point out that the peculiarity of our approach is given by the
possibility of looking at the conditional probability $P(A|B)$
as a primitive concept,with no need of assuming that the probability of
the {\em conditioning} event $B$ be positive.
On the contrary, in the probabilistic approaches usually adopted
in the literature, see, e.g., (Adams 1975, Hawthorne 1996, Schurz 1998),
by definition the quantity $P(A|B)$ is the {\em ratio} of $P(AB)$ and $P(B)$
if $P(B) \neq 0$, with $P(A|B) = 1$ if $P(B) = 0$. Notice that a clear
{\em rationale} of this latter assumption does not seem to exist.
Indeed, in the framework of coherence, this
assumption is not made and the case of {\em conditioning} events of zero
probability is managed without any problem: as an example, the
condition $P(A|B) + P(A^c|B) = 1$, where $A^c$ denotes the negation of $A$,
is satisfied also when $P(B) = 0$. We think that the opportunity, offered
by the probabilistic approach based on coherence, of developing a  completely
general treatment of probabilistic default reasoning is important specially in
the field of nonmonotonic reasoning where infinitesimal probabilities
play a significant role. Moreover, exploiting the algorithms developed in our framework,
the lower and upper probability bounds associated with the conditional assertions
of a given knowledge base can be propagated to further conditional assertions,
obtaining in all cases the tightest probability bounds.
Beside allowing a more flexible and realistic approach to probabilistic default reasoning,
this provides an exact illustration of the degradation of System P rules
when interpreted in probabilistic terms.
The probabilistic approach based on coherence has been adopted in  many
recent papers, see, e. g., (Biazzo, and Gilio 1999), (Capotorti, and
Vantaggi 1999), (Coletti 1994), (Coletti, and Scozzafava 1996),
(Coletti, and Scozzafava 1999), (Gilio 1995a), (Gilio 1995b), (Gilio 1999),
(Gilio, and Ingrassia 1998), (Gilio, and Scozzafava 1994), (Lad 1999),
(Lad, Dickey, and Rahman 1990), (Scozzafava 1994).
The algorithms described in (Biazzo, and Gilio 1999), based on the linear programming technique,
have been implemented with Maple V.
The paper is organized as follows. In Section 2 we give some preliminary concepts on
coherence and probability logic. In Section 3 we consider the definitions
of probabilistic consistency and entailment given by Adams; then we recall the main
inference rules of his probability logic. In Section 4 we consider the propagation of lower
and upper probability bounds in System P. We also examine the disjunctive Weak Rational
Monotony of System P$^+$. In Section 5 we examine the propagation of lower bounds with
real $\epsilon-$values. Then, in Section 6 we apply the results to (a slightly modified
version of) an example given in (Bourne, and Parsons 1998). Finally, in Section 7 we give some
conclusions.

\section{Some preliminaries}

We recall some preliminary concepts on coherence of imprecise probability assessments and
on probability logic. Given a family $\F_n = \{E_1 | H_1, \ldots, E_n|H_n \}$
and a vector $\A_n = (\alpha_1, \ldots, \alpha_n)$ of lower bounds
$P(E_i|H_i) \geq \alpha_i$, with $i \in J_n = \{1, \ldots, n\}$,
we consider the following definition of generalized coherence (g-coherence),
given in (Biazzo, and Gilio 1999), which essentially coincides with a previous one
introduced in (Gilio 1995a).
\begin{Def} \label{IMPR}{\rm
The vector of lower bounds $\A_n$ on $\F_n$ is said
g-coherent  if and only if there exists a precise coherent assessment
$\P_n = (p_1, \ldots, p_n)$
on $\F_n$, with $p_i = P(E_i|H_i)$, which is consistent with $\A_n$, that is
such that $p_i  \geq  \alpha_i$  for each $i \in J_n$.
}\end{Def}
The Definition ~\ref{IMPR} can be also applied to
imprecise assessments like
\[
\alpha_i \leq P(E_i|H_i) \leq  \beta_i \; , \; \; i \in J_n,
\]
since each inequality $P(E_i|H_i) \leq  \beta_i$ amounts to the inequality
$P(E_i^c|H_i) \geq  1-\beta_i$. \\
Then, given an imprecise assessment $\A_n= (\alpha_1, \ldots, \alpha_n)$
on $\F_n$, a suitable procedure, given in  (Gilio 1995b),
can be used to check the g-coherence of $\A_n$.
The g-coherent extension of $\A_n$ to a further conditional event
$E_{n+1}|H_{n+1}$ has been studied in (Biazzo, and Gilio 1999) where, defining a
suitable interval $[p_{\circ},p^{\circ}] \subseteq [0,1]$, the following
result has been obtained.
\begin{Thm}\label{GENdF}{\rm
Given a g-coherent imprecise assessment
$\A_n  = ([\alpha_i, \beta_i], i \in J_n)$ on the family
$\F_n = \{E_i|H_i, \; i \in J_n\}$, the extension
$[\alpha_{n+1}, \beta_{n+1}]$  of $\A_n$ to a further conditional event
$E_{n+1}|H_{n+1}$ is g-coherent if and only if the following condition is
satisfied
\begin{equation}\label{INTERS}
[\alpha_{n+1}, \beta_{n+1}] \cap [p_{\circ},p^{\circ}] \neq \emptyset \; .
\end{equation}
} \end{Thm}
In the quoted paper an algorithm has been proposed to determine
$[p_{\circ},p^{\circ}]$. Moreover, starting with a g-coherent assessment $\A_n$, by the same
algorithm it is possible to determine the corresponding assessment $\A^*_n$
coherent wrt definition given in (Walley 1991). \\
We can frame our approach to the problem of propagating imprecise
conditional probability assessments ({\em probabilistic deduction})
from the {\em probability logic} point of view, see, e. g.,
(Frisch, and Haddawy 1994), (Lukasiewicz 1998), (Nilsson 1986). \\
We associate to each conditional assertion  $H \; \s~ E$ in the
knowledge base $\K$ a probability interval $[\alpha,\beta]$.
In particular, given a family $\F_n$ of $n$ conditional assertions, we
consider an interval-valued probability
assessment $\A_n = ([\alpha_i,\beta_i], i = 1, \ldots, n)$.
Then, we can look at the pair $(\F_n, \A_n)$ as
a probabilistic knowledge base, where each imprecise assessment
$\alpha_i \leq P(E_i|H_i) \leq \beta_i$ is a probabilistic formula
denoted by $(E_i|H_i)[\alpha_i,\beta_i]$. In our approach a probabilistic
interpretation is just a coherent precise conditional probability assessment $\P_n$  on
$\F_n$. A probabilistic interpretation $\P_n = (p_1, \ldots, p_n)$ is a {\em
model} of a probabilistic formula $(E_i|H_i)[\alpha_i, \beta_i]$ iff ${\P}_n
\models (E_i|H_i)[\alpha_i, \beta_i]$, that is $\alpha_i \leq p_i \leq
\beta_i$. $\P_n$ is a model of the probabilistic knowledge base $\K = (\F_n,
\A_n)$, denoted $\P_n \models \K$, iff $\P_n \models (E|H)[\alpha, \beta]$ for
every $(E|H)[\alpha, \beta] \in \K$. Therefore, $\P_n$ is a model of $\K =
(\F_n, \A_n)$ iff $\P_n$ is consistent with $\A_n$. A set of probabilistic
formulas $\K$ is {\em satisfiable} iff a model of $\K$ exists, therefore the
concept of satisfiability of $\K = (\F_n, \A_n)$  coincides with that of
g-coherence of $\A_n$ on $\F_n$. A probabilistic formula
$(E_{n+1}|H_{n+1})[\alpha_{n+1}, \beta_{n+1}]$ is a {\em logical consequence}
of $\K = (\F_n, \A_n)$, denoted $\K \models (E_{n+1}|H_{n+1})[\alpha_{n+1},
\beta_{n+1}]$, iff
\[
\alpha_{n+1} \leq inf \; \I \; , \; \; \; \beta_{n+1} \geq sup \; \I \; ,
\]
where $\I$ is the set of the real values $p$ such that there exists a model of
$\K \cup \{(E_{n+1}|H_{n+1})[p, p]\}$. As shown by the condition
\(-~\ref{INTERS}), in our approach this amounts to
\[
[p_{\circ},p^{\circ}] \subseteq [\alpha_{n+1}, \beta_{n+1}].
\]
A probabilistic formula
$(E_{n+1}|H_{n+1})[\alpha_{n+1}, \beta_{n+1}]$ is a
{\em tight logical consequence} of $\K = (\F_n, \A_n)$, denoted
$\K \models_{tight} (E_{n+1}|H_{n+1})[\alpha_{n+1}, \beta_{n+1}]$,
iff
\[
\alpha_{n+1} = inf \; \I \; , \; \; \; \beta_{n+1} = sup \; \I \; ,
\]
that is
\[
\alpha_{n+1} = p_{\circ} \; , \; \; \; \beta_{n+1} = p^{\circ} \; .
\]
Considering a {\em probabilistic }{\em query} $(E_{n+1}|H_{n+1})[\alpha,
\beta]$, where $\alpha$ and $\beta$ are two different variables, to a
probabilistic knowledge base $\K = (\F_n, \A_n)$ a {\em correct answer} is any
$[\alpha, \beta] = [\alpha_{n+1}, \beta_{n+1}] \supseteq
[p_{\circ},p^{\circ}]$, that is such that $\K \models
(E_{n+1}|H_{n+1})[\alpha_{n+1}, \beta_{n+1}]$. The {\em tight answer} is
$[\alpha, \beta] = [p_{\circ},p^{\circ}]$.

\section{Probabilistic consistency and entailment}

We recall that in (Adams 1975) the
conditional assertion {\em "if} $A$ {\em then} $B$ {\em "} is looked at as
$P(B|A) \geq 1 - \epsilon \; \; \;  (\forall \epsilon > 0) $.
Adopting a more realistic point of view we may look at the same conditional assertion as the
probabilistic formula $(B|A)[\alpha, \beta]$, with $0 \leq \alpha \leq \beta \leq 1$,
where usually $\beta = 1$.
Then, a (probabilistic) knowledge base might be defined as a family of
probabilistic formulas $\K = \{(E|H)[\alpha, \beta]\}$.
In (Adams 1975) the following definition has been given.
\begin{Def}{\rm
The knowledge base $\cal K$ is {\em probabilistically
consistent} ({\em p-consistent}) if, for every $\epsilon > 0$,
there exists a  probability  $P$ on {\cal A}, proper for $\cal K$,
such that $P(E|H) \geq 1-\epsilon$ for every  $E|H  \in \cal K$ .
} \end{Def}
In our framework the concept of probabilistic consistency can be defined
in the following way.

\begin{Def}{\rm
The knowledge base $\cal K$ is {\em probabilistically
consistent} ({\em p-consistent}) if, for every set of lower bounds
$\A = \{\alpha_{E|H}, E|H \in \K\}$ on $\K$, there exists a precise coherent
probability assessment $P = \{p_{E|H}, E|H \in \K\}$ on $\K$, with $p_{E|H} = P(E|H)$,
such that, for each $E|H \in \K$, $p_{E|H} \geq  \alpha_{E|H}$.
} \end{Def}

We recall the concept of probabilistic entailment as defined in
(Adams 1975).
\begin{Def}{\rm
A p-consistent knowledge  base  $\cal K$  {\em probabilistically  entails}
({\em p-entails}) the conditional $B|A$ if, for every $\epsilon > 0$,
 there  exists $\delta > 0$ such that for all probabilities $P$,
proper for ${\cal K} \cup \{B|A\}$,  if $P(E|H) \geq 1 - \delta$
for each $E|H \in \cal K$, then $P(B|A) \geq 1 - \epsilon$.
} \end{Def}
In our framework the concept of probabilistic entailment can be defined in
the following way.

\begin{Def}\label{PE}{\rm
A p-consistent knowledge  base  $\cal K$  {\em probabilistically  entails}
the conditional $B|A$ if there exists a subfamily $\F \subseteq \K$ such
that, for every $\alpha_{B|A}$, there  exists a set of lower bounds
$\A = \{\alpha_{E|H}, E|H \in \F\}$ on $\F$ such that for all precise coherent
probability assessment $P = \{p_{B|A}, p_{E|H}, E|H \in \F\}$ on $\F \cup \{B|A\}$,
with $p_{B|A} = P(B|A), p_{E|H} = P(E|H)$, if $p_{E|H} \geq \alpha_{E|H}$
for each $E|H \in \F$, then $p_{B|A} \geq  \alpha_{B|A}$.
} \end{Def}
The probabilistic entailment of $B|A$ by $\cal K$ is denoted by the symbol
${\cal K} \Rightarrow B|A$.
In (Adams 1975) a  suitable  set  ${\cal R}$  of  seven  inference  rules
has  been introduced, by  means  of  which  an  inferential  system
can  be developed  to  deduce  in  an  automatic  way  all  the
plausible conclusions of a knowledge base $\cal K$.See also (Pearl 1988).
The fundamental rules  of the set $\cal R$ are the following ones:
\[
\begin{array}{ll}
R1. &  A \;  \s~ C \; , \; A \;  \s~ B \;  \Rightarrow \;
AB \;  \s~ C \\ & (Triangularity) \\ \ \\
R2. &  AB \; \s~ C \; , \; A \; \s~ B  \;  \Rightarrow \;
A \;  \s~ C \\ & (Bayes) \\ \ \\
R3 &  A  \; \s~ C \; , \; B \;  \s~ C \;  \Rightarrow \;
A \vee B  \; \s~ C \\ & (Disjunction)
\end{array}
\]
The previous rules, among others, have been used  (with  different
names and in the framework of symbolic  approaches  too)  by  many
authors, with the aim of developing some  nonmonotonic  logics  to
formalize the plausible reasoning (see e.g. (Kraus, Lehmann, and Magidor 1990),
(Lehmann, and Magidor 1992), (Dubois, and Prade 1994)).
See also the survey given in  (Benferhat, Dubois, and Prade 1997). \\
We recall that in (Gilio 1997) the following concept of
{\em strict probabilistic consistency} has been introduced.
\begin{Def}{\rm
The knowledge base $\cal K$ is strictly p-consistent if
the probability assessment $P$ on $\cal K$, such that $P(E|H) = 1$
for each $E|H \in \cal K$, is coherent.
}
\end{Def}
Then the following result, which has some relations with the Theorem 3
given in (Schurz 1998), has been given (without proof).
\begin{Thm} \label{STRICT} {\rm
 $\cal K$  is  p-consistent  if  and  only  if  $\cal K$
is  strictly p-consistent.
}
\end{Thm}
The proof of Theorem ~\ref{STRICT} is based on the following
three assertions:
\begin{description}
\item[a.]If $\cal K$ is strictly p-consistent, then $\cal K$ is
p-consistent;
\item[b.] If $\cal K$ is p-consistent, then $\cal K$ is consistent;
\item[c.] If $\cal K$ is consistent, then $\cal K$ is strictly
p-consistent.
\end{description}
Hence, the property of p-consistency can be  simply defined  on
the basis of the property of  strict  p-consistency. Moreover,  the
following well known result
\begin{Thm}{\rm
 If $\cal K$ is consistent, then ${\cal K} \Rightarrow B|A$
 if and  only  if ${\cal K} \cup \{B^c|A\}$ is inconsistent.
}
\end{Thm}
can  be reformulated in the following way:
\begin{Thm}{\rm
 Given a consistent knowledge base $\cal K$ and  a
conditional $B|A$, $\cal K$ p-entails $B|A$ if and only if the
probability  assessment  $P$ on ${\cal K} \cup \{B^c|A\}$, with
$P(B^c|A) = P(E|H) = 1$ for each $E|H \in \cal K$, is not coherent.
}
\end{Thm}

\section{Exact propagation of probability bounds in System P}
We recall that the inference rules in System P are the following ones:
\[
\begin{array}{ll}
1. & A \; \s~ A  \\
& (Reflexivity) \\ \ \\

2. & \models A \leftrightarrow B,  \; A \; \s~ C
\; \Longrightarrow \; B \; \s~ C  \\
 & (Left \; Logical \; Equivalence) \\ \ \\

3. & \models B \; \rightarrow C, \; A \; \s~ B   \;
\Longrightarrow \; A \; \s~ C  \\
 & (Right \; Weakening) \\ \ \\
4. & A \; \s~ B, \; A \; \s~ C \;
\Longrightarrow \; A \; \s~ B   C \\
 & (And) \\ \ \\
5. & A \; \s~ C, \; A \; \s~ B   \;  \Longrightarrow \;
     AB \; \s~ C \\
 & (Cautious \; Monotonicity) \\ \ \\
6. &
A \; \s~ C, \; B \; \s~ C   \;  \Longrightarrow \;
A \vee B \; \s~ C \\
 & (Or) \\
\end{array}
\]
Two derived rules in System P are
\[
\begin{array}{lll}
7. \hspace{0.5 cm}
A   B \; \s~ C, \; \; A \; \s~ B \;
\Longrightarrow \; A \; \s~ C & & ({\em Cut}) \\ \ \\
8. \hspace{0.5 cm}
A   B \; \s~ C \; \Longrightarrow \; A \; \s~ B \rightarrow C
& & ({\em S})
\end{array}
\]
As we can see, the rules $R1, R2, R3$ coincide respectively with the
rules {\em Cautious Monotonicity, Cut, Or} in System P.
In (Adam86) an extended probability logic (System P$^+$) was developed
by allowing the disjunction of conditionals in the conclusion of inferences
and by adding the following dWRM rule.
\[
\begin{array}{ll}
9. & A \; \s~ C \; \Longrightarrow
\; A \; \s~ B^c \; \; \vee \; \; AB \; \s~ C \\
 & (disjunctive \; Weak \; Rational \; Monotony)
\end{array}
\]
Now we will show how the probability
intervals associated with the antecedents in each rule of System P propagate
in a precise way to the consequent of the given rule. We will also examine the
rules Cut, S and dWRM. These bounds will provide
an exact illustration of the degradation of System P rules when interpreted in
probabilistic terms. Perhaps, as a by-product, also a better
appreciation of the results given in (Hawthorne 1996) on the interpretation of
nonmonotonic conditionals in terms of probabilistic support levels could be obtained.
We assume that the basic events $A, B, C$ are logically independent.
\begin{enumerate}
\item {\em Reflexivity rule}. As for every assertion $A$ the
assessment $P(A|A) = p$ is coherent if
and only if $p = 1$, to the conditional assertion $A \; \s~ A$ we associate
the interval  $[\alpha, \beta] = [1,1]$.
\item {\em Left Logical  Equivalence  rule}. If two assertions $A, B$ are
equivalent, then for every assertion $C$ the assessment $(x,y)$ on
the pair of conditional events $C|A, C|B$ is coherent if and only $x=y$.
Therefore we associate the same probability interval to the conditional
assertions $C|A, C|B$. In other words, the assessment $[\alpha, \beta]$ on
$C|A$ propagates to the same interval on $C|B$.
\item {\em Right Weakening rule}. If $B \subseteq C$ then, defining the inclusion
among conditional events as in (Goodman, and Nguyen 1988),
it is $B|A \subseteq C|A$ and then the assessment $(x,y)$ on
the pair of conditional events $B|A, C|A$ is coherent if and only
$x \leq y$, see (Gilio 1993). Therefore, the assessment $[\alpha, \beta]$
on $B|A$ propagates to the interval $[\alpha, 1]$ on $C|A$.
\item  {\em And rule}. Given the assessment
$(x,y)$ on the pair of conditional events $B|A, C|A$, as well known
the extension $P(BC|A) = z$ is coherent if and only if
\[
Max \; \{0, x + y - 1\} = z' \leq z \leq z'' = Min \; \{x,y\}.
 \]
Therefore, the probability intervals
$[\alpha_1, \beta_1], [\alpha_2, \beta_2]$ on the antecedents $B|A, C|A$
of the rule propagate to the exact interval $[\alpha_3, \beta_3]$, with
\begin{equation}\label{AND}
\begin{array}{l}
\alpha_3 = Min \; z' = Max \; \{0, \alpha_1 + \alpha_2 - 1\} \; , \\ \ \\
\beta_3 = Max \; z'' = Min \; \{\beta_1, \beta_2\} \; ,
\end{array}
\end{equation}
on the consequent $BC|A$.
\item {\em Cautious Monotonicity rule}. Given the assessment
$(x,y)$ on the pair of conditional events $C|A, B|A$, as proved in (Gilio 1995b),
the extension $P(C|AB) = z$ is coherent if and only if $z \in [z', z'']$, with

\[
\begin{array}{l}
z' \; = \; \left\{\begin{array}{lcl}
\frac{x+y-1}{y} \; , & \mbox{if} & x+y > 1  \\
0 \; , & \mbox{if} & x+y \leq 1
\end{array} \right. \; ,  \\ \ \\
z'' \; = \; \left\{\begin{array}{lcl}
\frac{x}{y} \; , & \mbox{if} & x < y  \\
1 \; , & \mbox{if} & x \geq y
\end{array} \right. \; .
\end{array}
\]
We observe that, for $(x,y) \in [\alpha_1, \beta_1] \times [\alpha_2, \beta_2]$
the function $f(x,y)$ attains its minimum value at the point $(\alpha_1, \alpha_2)$.
Then, it follows
\begin{equation}
\begin{array}{ll}
\alpha_3  \: = \:
\left\{\begin{array}{lll}
\frac{\alpha_1+\alpha_2-1}{\alpha_2},  & \mbox{if } &
\alpha_1 + \alpha_2 > 1  \\
0,  & \mbox{if } &  \alpha_1 + \alpha_2 \leq 1
\end{array} \right. \label{ACM}
\end{array}
\end{equation}
Moreover, the function $g(x,y)$ attains its maximum
value at the point $(\beta_1, \alpha_2)$. Then it follows
\begin{equation}
\beta_3  \: = \: \left\{\begin{array}{lll}
\frac{\beta_1}{\alpha_2},  & \mbox{if } & \beta_1 < \alpha_2 \\
1,  & \mbox{if } &  \beta_1 \geq \alpha_2
\end{array} \right.  \label{zeta^*}
\end{equation}
Then, in order to determine the interval $[\alpha_3, \beta_3]$ we have to
consider the position of the vertices $(\alpha_1, \alpha_2), (\beta_1, \alpha_2)$ wrt
diagonals of the unitary square $[0,1]^2$.
\item {\em Or rule}. Given the assessment
$(x,y)$ on the pair of conditional events $C|A, C|B$, it can be proved that
the extension $P(C|(A \vee B)) = z$ is
coherent if and only if $z \in [z', z'']$, with
\[
\begin{array}{l}
z' \; = \; \left\{\begin{array}{lll}
\frac{xy}{x+y-xy}, & \mbox{if} & (x,y) \neq (0,0) \\
0 , & \mbox{if} & (x,y) = (0,0)
\end{array}
\right. \; ,  \\ \ \\
z'' \; = \; \left\{\begin{array}{lll}
\frac{x+y-2xy}{1- xy}, & \mbox{if} & (x,y) \neq (1,1) \\
1 ,  & \mbox{if} & (x,y) = (1,1)
\end{array}
\right. \; .
\end{array}
\]
Moreover, we observe that both $z'$ and $z''$ increase as either $x$ or $y$
increase. Therefore, the probability intervals
$[\alpha_1, \beta_1], [\alpha_2, \beta_2]$ on the antecedents $C|A, C|B$
of the rule propagate, under the condition
$(\alpha_1,\alpha_2) \neq (0,0), \; (\beta_1,\beta_2) \neq (1,1)$,
to $[\alpha_3, \beta_3]$, with
\begin{equation}\label{OR}
\alpha_3 \; = \; \frac{\alpha_1\alpha_2}{\alpha_1+\alpha_2-\alpha_1\alpha_2} \; ,
\end{equation}
\begin{equation}\label{ORR}
\beta_3 \; = \; \frac{\beta_1+\beta_2-2\beta_1\beta_2}{1- \beta_1\beta_2} \; ,
\end{equation}
on the consequent $C|(A \vee B)$.

\end{enumerate}

Concerning the rules Cut and S we have the following results. \\
(e) {\em Cut rule}. Given the assessment
$(x,y)$ on the pair of conditional events $C|AB, B|A$, it can be proved that
the extension $P(C|A) = z$ is coherent if and only if
\[
xy \leq z \leq xy + 1 - y \; .
\]
Therefore, the probability intervals
$[\alpha_1, \beta_1], [\alpha_2, \beta_2]$ on the antecedents $C|AB, B|A$
of the rule propagate to $[\alpha_3, \beta_3]$, with
\begin{equation}\label{CUT}
\alpha_3  =
\alpha_1\alpha_2 \; , \; \; \;
\beta_3 = \beta_1\alpha_2 + 1 - \alpha_2 \; ,
\end{equation}
on the consequent $C|A$. \\
(f) {\em S rule}. As $C|AB \subseteq (B^c \vee C)|A$, then the assessment
$(x,y)$ on the conditional events $C|AB, (B^c \vee C)|A$ is coherent if and
only if $x \leq y$.
Therefore, the probability interval $[\alpha_1, \beta_1]$ on the antecedent
$C|AB$ of the rule propagates to $[\alpha_2, \beta_2]$, with
\[
\alpha_2  = \alpha_1, \; \; \; \beta_2 = 1 \; ,
\]
on the consequent $(B^c \vee C)|A$. \\
(g) {\em dWRM rule}. Let $\P = (x,y,z)$ a probability assessment on
the family $\{C|A, B^c|A, C|AB\}$. For this family the constituents
(possible worlds) are
\[
\begin{array}{l}
C_0 = A^c \; ,  \; \; C_1 = ABC \; , \; \; C_2 = ABC^c \; , \\ \ \\
C_3 = AB^cC \; , \; \; C_4 = AB^cC^c \; .
\end{array}
\]
To the constituents $C_1, \ldots, C_4$ we associate the points
\[
\begin{array}{l}
Q_1 = (1,0,1) \; , \; \; Q_2 = (0,0,0) \; , \\ \ \\
Q_3 = (1,1,z) \; , \; \;  Q_4 = (0,1,z) \; .

\end{array}
\]
Then, based on the method given in (Gilio 1995b) and denoting by $\I$
the convex hull of the points $Q_1, \ldots, Q_4$, it can be proved that
the coherence of $\P$ amounts to the condition $\P \in \I$. Notice that
in general this condition is necessary but not sufficient for the coherence
of an assessment $\P_n = (p_1, \ldots, p_n)$ on a family
$\F_n = \{E_1|H_1, \ldots, E_n|H_n\}$. The study of the condition
$\P \in \I$ requires considering the equations of the four planes determined
respectively by the terns of points
\[
\begin{array}{l}
\{Q_1, Q_2, Q_3\} \; , \; \; \; \{Q_2, Q_3, Q_4\} \; , \\ \ \\
\{Q_1, Q_2, Q_4\} \; , \; \; \; \{Q_1, Q_3, Q_4\} \; .
\end{array}
\]
Denoting by $X,Y,Z$ the axes' coordinates, the equations are given
respectively by
\[
\begin{array}{l}
Z = X + (z-1)Y \; , \; \; \; \; Z = zY \; , \\ \ \\
Z = X + zY \; , \; \; \; \; Z = (z-1)Y + 1 \; .
\end{array}
\]
Then, given the values $x,y$, it is
\[
\P \in \I \; \Longleftrightarrow \; z' \leq z \leq z'' \; ,
\]
where
\[
\begin{array}{l}
z' \; = \;  \left\{\begin{array}{lcl}
\frac{x-y}{1-y} \; , & \mbox{if} & x > y  \\
0 \; , & \mbox{if} & x \leq y
\end{array} \right.  \; , \\ \ \\
z'' \; = \; \left\{\begin{array}{lcl}
\frac{x}{1-y} \; , & \mbox{if} & x+y < 1  \\
1 \; , & \mbox{if} & x+y \geq 1
\end{array} \right. \; .
\end{array}
\]
In order to examine the probabilistic interpretation of the rule
we introduce a partition $\{{\R}_1, {\R}_2, {\R}_3, {\R}_4\}$ of the unitary
square $[0,1]^2$, with
\[
\begin{array}{lll}
{\R}_1 & = & \{(x,y): x+y < 1, x \geq y \} \; , \\ \\
{\R}_2 & = & \{(x,y): x+y < 1, x < y \} \; , \\ \ \\
{\R}_3 & = & \{(x,y): x+y \geq 1, x < y \} \; , \\ \ \\
{\R}_4 & = & \{(x,y): x+y \geq 1, x \geq y \} \; .
\end{array}
\]
We have to examine the case in which $x$ is "high", therefore ${\R}_2$
is not of interest. In ${\R}_3$, since $x < y$, if $x$ is "high" then $y$ is
"high" too. In ${\R}_1$ and ${\R}_4$ it is $z \geq z' = \frac{x-y}{1-y}$
so that, if $x$ is "high" and $y$ is "not high", then  $z$ is "high". \\
Concerning propagation of probability intervals, if we consider
the assessments $[\alpha_1, \beta_1], [\alpha_2, \beta_2]$ on the
conditional events $C|A, B^c|A$, then for the interval $[\alpha_3, \beta_3]$
associated with $C|AB$ we first observe that the quantity $\frac{x}{1-y}$
attains its maximum value at the point $(\beta_1, \alpha_2)$, while the
quantity $\frac{x-y}{1-y}$ attains its minimum value at the point
$(\alpha_1, \beta_2)$. Then, we have:
\begin{equation}
\alpha_3 \; = \;
\left\{\begin{array}{lll}
\frac{\alpha_1-\beta_2}{1-\beta_2}, & \mbox{if} & \alpha_1 \geq \beta_2 \\
0 , & \mbox{if} & \alpha_1 < \beta_2
\end{array} \right.
\end{equation}
\begin{equation}
\beta_3 \; = \;
\left\{\begin{array}{lll}
\frac{\beta_1}{1-\alpha_2}, & \mbox{if} & \beta_1 + \alpha_2 < 1  \\
1 , & \mbox{if} & \beta_1 + \alpha_2 \geq 1
\end{array} \right.
\end{equation}

\begin{Rem}{\rm
Using the lower bounds computed previously, we can verify the probabilistic entailment
in each inference rule on the basis of Definition ~\ref{PE}. We have
\begin{itemize}
\item {\em And rule.} For each given value $\alpha_3$, from \(-~\ref{AND}) we have
that, for every $(\alpha_1, \alpha_2) \in [\alpha_3,1] \times [\alpha_3,1]$ such that
$\alpha_1 + \alpha_2 = 1 + \alpha_3$, if $P(B|A) \geq \alpha_1, P(C|A) \geq \alpha_2$
then $P(BC|A) \geq \alpha_3$.
\item {\em Cautious Monotonicity rule.}
For each given value $\alpha_3$, from \(-~\ref{ACM}) we have
that, for every $(\alpha_1, \alpha_2) \in [\alpha_3,1] \times (0,1]$ such that
$\alpha_1 + (1 - \alpha_3)\alpha_2 = 1$, if $P(C|A) \geq \alpha_1, P(B|A) \geq \alpha_2$
then $P(C|AB) \geq \alpha_3$.
\item {\em Or rule.} For each given value $\alpha_3$, from \(-~\ref{OR}) we have
that, for every $(\alpha_1, \alpha_2) \in [\alpha_3,1] \times [\alpha_3,1]$ such that
$\alpha_2 = \frac{\alpha_1\alpha_3}{\alpha_1(1+\alpha_3)-\alpha_3}$,
if $P(C|A) \geq \alpha_1, P(C|B) \geq \alpha_2$ then $P(C|A \vee B) \geq \alpha_3$.
\item {\em Cut rule.} For each given value $\alpha_3$, from \(-~\ref{CUT}) we have
that for every $(\alpha_1, \alpha_2) \in [\alpha_3,1] \times [\alpha_3,1]$ such that
$\alpha_2 = \frac{\alpha_3}{\alpha_1}$, if $P(C|AB) \geq \alpha_1, P(B|A) \geq \alpha_2$
then $P(C|A) \geq \alpha_3$.
\end{itemize}
}\end{Rem}

\section{Propagation with $\epsilon-$values}
The results of the previous section can be examined in the particular case
in which for $i = 1,2$ it is $[\alpha_i, \beta_i] = [1 - \epsilon_i, 1]$.
As it can be verified, the $\epsilon-$values propagate in the following way.
\begin{itemize}
\item {\em And rule}. From \(-~\ref{AND}), the probability bounds \linebreak
$[1 - \epsilon_1, 1], [1 - \epsilon_2, 1]$ on the antecedents $B|A, C|A$
of the rule propagate, on the consequent $BC|A$, to the exact bounds
$[1 - \epsilon_3, 1]$, with
\begin{equation}\label{AND2}
\epsilon_3 = \epsilon_1 + \epsilon_2
\end{equation}

\item {\em Cautious Monotonicity rule}.
From \(-~\ref{ACM}), the probability intervals $[1 - \epsilon_1, 1]$,
$[1 - \epsilon_2, 1]$ on the antecedents $C|A, B|A$
of the rule propagate, on the consequent $C|AB$,
to $[1 - \epsilon_3, 1]$, with
\begin{equation}\label{CM}
\epsilon_3 =  \frac{\epsilon_1}{1 - \epsilon_2}
\end{equation}

\item {\em Or rule}.
From \(-~\ref{OR}), the probability intervals  \linebreak
$[1 - \epsilon_1, 1], [1 - \epsilon_2, 1]$ on the antecedents $C|A, C|B$
of the rule propagate, on the consequent $C|(A \vee B)$,
to  \linebreak $[1 - \epsilon_3, 1]$, with
\begin{equation}\label{OR2}
\epsilon_3 =  \frac{\epsilon_1 + \epsilon_2 - 2\epsilon_1\epsilon_2}
{1 - \epsilon_1\epsilon_2}
\end{equation}

\item {\em Cut rule}.
From \(-~\ref{CUT}), the probability intervals \linebreak
$[1 - \epsilon_1, 1], [1 - \epsilon_2, 1]$ on the antecedents $C|AB, B|A$
of the rule propagate, on the consequent $C|A$,
to  \linebreak $[1 - \epsilon_3, 1]$, with
\begin{equation}\label{CUT2}
\epsilon_3 =  \epsilon_1 + \epsilon_2 - \epsilon_1\epsilon_2
\end{equation}\end{itemize}
\begin{Rem}{\rm
Our results concerning the value of $\epsilon_3$ coincide with that
ones obtained in (Bourne, and Parsons 1998) for the rules {\em And} and
{\em Cautious Monotonicity} and are better for the rules {\em Or} and
{\em Cut}, as from \(-~\ref{OR2}) and \linebreak \(-~\ref{CUT2}) one has respectively
\[
\begin{array}{cc}
\epsilon_3 \; = \; \frac{\epsilon_1 + \epsilon_2 - 2\epsilon_1\epsilon_2}
{1 - \epsilon_1\epsilon_2}  \; < \;  \epsilon_1 + \epsilon_2  &
(\epsilon_1 < 1, \epsilon_2 < 1) , \\ \ \\
\epsilon_3  \; = \;  \epsilon_1 + \epsilon_2 - \epsilon_1\epsilon_2  \; < \;
\epsilon_1 + \epsilon_2  &  (\epsilon_1 > 0, \epsilon_2 > 0) .
\end{array}
\]
The use of the precise bounds may have some relevance when the inference rules
are applied with real $\epsilon-$values.
}
\end{Rem}

\section{An application}
We will now examine an example to give an idea, on one hand, of how much the
conclusions may be sensible
to the use of methods of exact propagation of probability bounds and, on another
hand, of the related phenomenon of degradation
of inference rules when interpreted in probabilistic terms.
The example is a modified version of an application
considered in (Bourne, and Parsons 1998) which was inspired by examples given in (Kraus, Lehmann, and Magidor 1990).
We consider a probabilistic knowledge base consisting of some conditional assertions,
which concern the fact that a given party has various attributes (the party is {\em great,
noisy}, {\em Linda} and {\em Steve} are {\em present}, and so on).
By the symbol
$A \; \s~_{\epsilon}B$ we denote the assessment $P(B|A) \geq 1 - \epsilon$.
We start with a knowledge
base which has the following rules and $\epsilon-$values:
\[
\begin{array}{ccc}
1. & & Linda \; \s~_{0.05} great \\ \ \\
2. & & Linda \; \s~_{0.2} Steve \\ \ \\
3. & & Linda \wedge Steve \; \s~_{0.1} \neg noisy \\ \ \\
4. & & Steve \; \s~_{0.05} Linda \\ \ \\
5. & & \neg noisy \; \s~_{0.2}  \neg great
\end{array}
\]
Notice that the conditional $"Linda \; \s~_{0.05} great"$ means that the
probability of the conditional event
\begin{center}
"(The party will be great $|$ Linda goes to the party)"
\end{center}
is greater than or equal to $1 - 0.05$, and so on.
We are interested in propagating the previous bounds to find the
$\epsilon-$values of the following conditionals:
\[
\begin{array}{ccc}\label{SET}
(a) & & Linda \; \s~_{\epsilon} \neg noisy \\ \ \\
(b) & & \top \; \s~_{\epsilon}\neg Linda \\ \ \\
(c) & & Linda \wedge Steve \; \s~_{\epsilon} great \wedge \neg noisy \\ \ \\
(d) & & Steve \; \s~_{\epsilon} \neg noisy \\ \ \\
(e) & & Linda \vee Steve \; \s~_{\epsilon} \neg noisy
\end{array}
\]
By the symbol $\top$ we denote (any tautology representing) the certain
event. \\ Applying the Cut rule to the conditionals
\[
Linda \; \s~_{0.2} Steve \; , \; \; \;
Linda \wedge Steve \; \s~_{0.1} \neg noisy \; ,
\]
we obtain the conditional
\[
Linda \; \s~_{0.28} \neg noisy \; .
\]
Then, applying the And rule  to the conditionals
\[
Linda \; \s~_{0.28} \neg noisy \; , \; \; \; Linda \; \s~_{0.05} great \; ,
\]
we obtain the conditional
\[
Linda \; \s~_{0.33} great \wedge \neg noisy \; .
\]
Applying the S rule to
\[
Linda \; \s~_{0.33} great \wedge \neg noisy
\]
we obtain
\[
\top \; \s~_{0.33} \neg Linda \vee great \wedge \neg noisy \; .
\]
Applying the S rule to
\[
\neg noisy \; \s~_{0.2}  \neg great
\]
we obtain
\[
\top \; \s~_{0.2} noisy \vee \neg great \; .
\]
Finally, applying the And rule to the conditionals
\[
\top  \; \s~_{0.33} \neg Linda \vee great \wedge \neg noisy , \;
\top  \; \s~_{0.2} noisy \vee \neg great ,
\]
we obtain the conditional
\[
\top  \; \s~_{0.53} \neg Linda \wedge (noisy \vee \neg great).
\]
Then, by the Right Weakening rule we have
\[
\top  \; \s~_{0.53} \neg Linda \wedge (noisy \vee \neg great)
\; \Longrightarrow \;
\top  \; \s~_{0.53} \neg Linda.
\]
Concerning the conditional $(c)$, applying the Cautious Monotonicity rule to
\[
Linda \; \s~_{0.05} great \; , \; \; \;
Linda \; \s~_{0.2} Steve \; ,
\]
we obtain
\[
Linda \wedge Steve \; \s~_{0.0625} great \; .
\]
Then, applying the And rule to the conditionals
\[
Linda \wedge Steve \; \s~_{0.0625} great, \;
Linda \wedge Steve \; \s~_{0.1} \neg noisy ,
\]
we obtain the conditional
\[
Linda \wedge Steve \; \s~_{0.0725} great \wedge \neg noisy \; .
\]
Concerning the conditional $(d)$, applying the Cut rule to the conditionals
\[
Linda \wedge Steve \; \s~_{0.1} \neg noisy \; , \; \;
Steve \; \s~_{0.05} Linda \; ,
\]
we obtain the conditional
\begin{equation}\label{15}
Steve \; \s~_{0.145} \neg noisy
\end{equation}
Then, applying the Or rule to the conditionals
\[
Steve  \s~_{0.145} \neg noisy \; , \; \;
Linda  \s~_{0.28} \neg noisy \; ,
\]
we obtain the conditional
\begin{equation}\label{16}
Linda \vee Steve  \s~_{0.358} \neg noisy
\end{equation}
We observe that, propagating the bounds with
$\epsilon_3 = \epsilon_1 + \epsilon_2$, instead of the
conditionals \(-~\ref{15}) and \(-~\ref{16}) we would obtain
respectively
\[
Steve \; \s~_{0.15} \neg noisy \; ,
\]
and
\[
Linda \vee Steve \; \s~_{0.425} \neg noisy \; .
\]

\section{Conclusions}
In this paper the inference rules of System P have been considered in the
framework of coherence.
We have also examined the {\em disjunctive Weak Rational Monotony}
proposed by Adams in his extended probability logic, corresponding to
System P$^+$. Differently from the probabilistic approaches generally given
in the literature, see, in particular, (Hawthorne 1996) and (Schurz 1998),
within our framework we can directly manage conditional probability
assessments, even if some (or possibly all the) conditioning events have
zero probability. We think that this opportunity is important specially in
the field of nonmonotonic reasoning where infinitesimal probabilities
play a significant role. Moreover, exploiting our algorithms, the lower
and upper probability bounds
associated with the conditional assertions of a given knowledge base
can be propagated to further conditional assertions,
obtaining in all cases the precise probability bounds. In particular,
beside allowing a more flexible
and realistic approach to probabilistic default reasoning,
this provides an exact illustration of the degradation of System P rules
when interpreted in probabilistic terms.   \\ \ \\
\subsection{Acknowledgments}
The author is grateful to the referees for the very helpful criticisms and suggestions.

\bibliographystyle{aaai}

\section{References}

Adams, E. W. 1975. {\em The Logic of  Conditionals}.
Dordrecht, Netherlands: Reidel. \\
Adams, E. W. 1986. On the logic of high probability,
{\em Journal of Philosophical Logic} 15: 255-279. \\
Benferhat, S.; Dubois, D.; and Prade, H. 1997. Nonmonotonic reasoning,
conditional objects and possibility theory. {\em Artificial Intelligence}
92(1-2): 259-276. \\
Biazzo, V.; and Gilio A. 1999. A generalization of the fundamental theorem of de Finetti for
imprecise conditional probability assessments.
{\em International Journal of Approximate Reasoning}. Forthcoming. \\
Bourne, R.; and Parsons, S. 1998. Propagating probabilities in System P.
In {\em Proceedings of the 11th International FLAIRS Conference}, 440-445. \\
Capotorti, A.; and Vantaggi, B. 1999. A general interpretation of
conditioning and its implication on coherence. {\em Soft Computing}
3/3: 148-153. \\
Coletti, G. 1994. Coherent   numerical   and   ordinal
probabilistic  assessments. {\em IEEE  Trans.  on  Systems,  Man,  and
Cybernetics} 24(12): 1747-1754. \\
Coletti, G.; and Scozzafava, R. 1996. Characterization   of
coherent conditional probabilities as a tool for their  assessment
and extension. {\em Journal   of   Uncertainty,    Fuzziness    and
Knowledge-based Systems} 4(2): 103-127. \\
Coletti, G.; and Scozzafava, R. 1999. Conditioning and inference in
intelligent systems. {\em Soft Computing} 3/3: 118-130. \\
Dubois, D.; and  Prade, H. 1991. Semantic consideration on order of
magnitude reasoning. In {\em Decision support systems and qualitative
reasoning}, Singh, M. G.; and Trave-Massuyes L. eds., 223-228.
Elsevier Science Publishers B. V., North Holland. \\
Dubois,  D.; and  Prade,  H. 1994.   Conditional   Objects   as
Nonmonotonic Consequence Relationships. {\em IEEE  Transactions  on  Systems,
Man, and Cybernetics}, 24(12): 1724-1740. \\
Frisch, A. M.; and Haddawy, P. 1994. Anytime Deduction for Probabilistic Logic,
{\em Artificial Intelligence} 69: 93-122. \\
Gilio, A. 1993. Conditional events and subjective probability in management
of uncertainty. In {\em Uncertainty  in  Intelligent  Systems - IPMU'92},
Bouchon-Meunier, B.; Valverde, L.; and Yager, R. R. eds., 109-120.
Elsevier Science Publ. B. V., North-Holland. \\
Gilio, A. 1995a. Probabilistic consistency of conditional probability
bounds. In {\em Advances in Intelligent Computing - IPMU '94,
Lecture Notes in Computer Science} 945, Bouchon-Meunier, B.; Yager, R. R.;
and Zadeh, L. A. eds., 200-209. Berlin Heidelberg: Springer-Verlag. \\
Gilio, A. 1995b. Algorithms for precise and imprecise conditional
probability assessments. In {\em Mathematical Models for Handling Partial
Knowledge in Artificial Intelligence}, Coletti, G.; Dubois,  D.; and
Scozzafava, R. eds., 231-254. New York: Plenum Press. \\
Gilio, A. 1997. Probabilistic modelling of uncertain conditionals. In
{\em Procedings of European Symposium on Intelligent Techniques}, 54-57.
Bari, Italy. \\
Gilio, A. 1999. Probabilistic relations among logically
dependent conditional events, {\em Soft Computing} 3/3: 154-161. \\
Gilio, A.; and Ingrassia, S. 1998. Totally coherent set-valued
probability assessments, {\em Kybernetika} 34(1): 3-15. \\
Gilio, A.; and Scozzafava, R. 1994. Conditional events in probability
assessments and revision, {\em IEEE Trans. on Systems, Man and
Cybernetics} 24(12): 1741-1746. \\
Goodman, I. R.; Nguyen H. T. 1988. Conditional objects and the
modeling of uncertainties. In {\em Fuzzy Computing Theory, Hardware and
Applications}, Gupta, M. M.; and Yamakawa, T. eds., 119-138.
New York: North-Holland. \\
Hawthorne, J. 1996. On the logic of nonmonotonic conditionals and
conditional probabilities, {\em Journal of Philosophical Logic} 25: 185-218. \\
Kraus, K.; Lehmann, D.; and  Magidor, M. 1990.  Nonmonotonic
reasoning, preferential models and cumulative  logics,  {\em Artificial
Intelligence}  44: 167-207. \\
Lad, F. 1999. Assessing the foundations for Bayesian networks:
a challenge to the principles and the practice, {\em Soft Computing}
3/3: 174-180. \\
Lad, F. R.; Dickey, J. M.; and Rahman, M. A. 1990. The fundamental theorem of
prevision, {\em Statistica}, anno L(1): 19-38. \\
Lehmann, D.; and Magidor, M. 1992.  What  does  a  conditional
knowledge base entail?, {\em Artificial Intelligence} 55: 1-60. \\
Lukasiewicz, T. 1998. Magic inference rules for probabilistic deduction
under taxonomic knowledge. In {\em Proceedings of the 14th Conference on
Uncertainty in  Artificial Intelligence}, 354-361. Madison, Wisconsin. \\
Nilsson, N. J. 1986. Probabilistic logic, {\em Artificial Intelligence} 28: 71-87. \\
Pearl, J. 1988.  {\em Probabilistic  Reasoning  in  Intelligent Systems: Networks
of Plausible Inference}. San Mateo, CA:   Morgan Kaufmann. \\
Raiman, O. (1989). Le raisonnement sur les ordres de grandeur, {\em Revue
d'Intelligence Artificielle} 3(4): 55-67.  \\
Schurz, G. 1998. Probabilistic semantics for Delgrande's conditional logic
and a counterexample to his default logic, {\em Artificial Intelligence}
102(1): 81-95. \\
Scozzafava, R. 1994. Subjective probability versus belief functions in
artificial intelligence, {\em Int. J. General Systems} 22: 197-206. \\
Walley, P. 1991. {\em Statistical reasoning with imprecise probabilities}.
London: Chapman and Hall.

\end{document}